\documentclass[12pt]{amsart}
\usepackage{a4}
\usepackage{amsthm, amsfonts, amssymb,latexsym}
\usepackage{enumerate, color}
\usepackage[english]{babel}
\usepackage[latin1]{inputenc}
\usepackage{comment }

\newtheorem{lemma}{Lemma}[section]

\newtheorem{theorem}{Theorem}[section]

\newtheorem{cor}{Corollary}[section]

\newcommand{\la}{\lambda}

\numberwithin{equation}{section}

\newcommand{\beq}[1]{\begin{equation}\label{#1}}
\newcommand{\eeq}{\end{equation}}

\title[An improved bound for the size of the set $A/A+A$]{An improved bound for the size of the set $A/A+A$}

\author[Oliver Roche-Newton]{Oliver Roche-Newton}
\address{O. Roche-Newton: RICAM, 69 Altenberger Stra{\ss}e, Johannes Kepler Universit\"{a}t, Linz, Austria }
\email{o.rochenewton@gmail.com }

\begin{document}

\begin{abstract}
It is established that for any finite set of positive real numbers $A$, we have
$$|A/A+A| \gg \frac{|A|^{\frac{3}{2}+\frac{1}{26}}}{\log^{1/2}|A|}.$$
\end{abstract} 
\maketitle
\section{Introduction}

Given a set $A$, we define its \textit{sum set}, \textit{product set} and \textit{ratio set} as
$$A+A:=\{a+b:a,b \in A\},\,\,\,\,\,\,\, AA:=\{ab:a,b \in A\},\,\,\,\,\,\,\, A/A:=\{a/b:a,b \in A, b \neq 0\}.$$
respectively. It was conjectured by Erd\H{o}s and Szemer\'{e}di that, for any finite set $A$ of integers, at least one of the sum set or product set has near-quadratic growth. Solymosi \cite{solymosi} used a beautiful and elementary geometric argument to prove that, for any finite set $A \subset \mathbb R$,
\begin{equation}
\max \{|A+A|,|AA|\} \gg \frac{|A|^{4/3}}{\log^{1/3}|A|}.
\label{soly}
\end{equation}
Recently, a breakthrough for this problem was achieved by Konyagin and Shkredov \cite{KS}. They adapted and refined the approach of Solymosi, whilst also utilising several other tools from additive combinatorics and discrete geometry in order to prove that
\begin{equation}
\max \{|A+A|,|AA|\} \gg |A|^{\frac{4}{3}+\frac{1}{20598}-o(1)}.
\label{KS}
\end{equation}
Further refinements in \cite{KS2} and more recently \cite{RSS} have improved this exponent to $4/3 +1/1509 +o(1)$. See \cite{KS}, \cite{KS2}, \cite{RSS} and the references contained therein for more background on the sum-product problem.

In this paper the related problem of establishing lower bounds for the sets
$$AA+A:=\{ab+c:a,b,c \in A\},\,\,\,\,\,\,\,\,\, A/A+A:=\{a/b+c:a,b,c \in A, b \neq 0\}$$
are considered. It is believed, in the spirit of the Erd\H{o}s-Szemer\'{e}di conjecture, that these sets are always large. It was conjectured by Balog \cite{balog} that, for any finite set $A$ of real numbers, $|AA+A| \geq |A|^2$. In the same paper, he proved the following result in that direction:

\begin{theorem} \label{thm:old}
Let $A$ and $B$ be finite sets of positive real numbers. Then
$$|AB+A| \gg |A||B|^{1/2}.$$
In particular,
$$|AA+A| \gg |A|^{3/2} \,\,\,\,\,\,\, |A/A+A| \gg |A|^{3/2}.$$
\end{theorem}

The proof of Theorem \ref{thm:old} uses similar elementary geometric arguments to those of \cite{solymosi}. In fact, one can obtain the same bound by a straightforward application of the Szemer\'{e}di-Trotter Theorem (see \cite[Exercise 8.3.3]{tv}). 

Some progress in this area was made by Shkredov \cite{shkredov}, who built on the approach of Balog in order to prove the following result: 

\begin{theorem} \label{thm:ilya}
For any finite set $A$ of positive real numbers,
\begin{equation}
|A/A+A| \gg  \frac{|A|^{\frac{3}{2}+\frac{1}{82}}}{\log^{\frac{2}{41}}|A|},
\label{a+a:a}
\end{equation}
\end{theorem}

The main result of this paper is the following improvement on Theorem \ref{thm:ilya}: 
\begin{theorem} \label{thm:main} Let $A$ be a finite set of positive reals. Then
$$|A/A+A| \gg \frac{|A|^{\frac{3}{2}+\frac{1}{26}}}{\log^{5/6}|A|}.$$
\end{theorem}

For the set $AA+A$ the situation is different, and it has proven rather difficult to beat the threshold exponent of $3/2$. A detailed study of this set can be found be in \cite{RRSS}. However, the corresponding problem for sets of integers is resolved, up to constant factors, thanks to a nice argument of George Shakan.\footnote{See http://mathoverflow.net/questions/168844/sum-and-product-estimate-over-integers-rationals-and-reals, where this argument first appeared.}

\subsection{Sketch of the proof of Theorem \ref{thm:main}} The proof is a refined version of the argument used by Balog to prove Theorem \ref{thm:old}. Balog's argument goes roughly as follows:

Consider the point set $A \times A$ in the plane. Cover this point set by lines through the origin. Let us assume for simplicity that all of these lines are equally rich, so we have $k$ lines with $|A|^2/k$ points on each line. Label the lines $l_1,l_2,\dots,l_k$ in increasing order of steepness. Note that if we take the vector sum of a point on $l_i$ with a point on $l_{i+1}$, we obtain a point which has slope in between those of $l_i$ and $l_{i+1}$. The aim is to show that many elements of $(A/A+A) \times (A/A+A)$ can be obtained by studying vector sums from neighbouring lines.

Indeed, for any $1 \leq i \leq k-1$, consider the sum set
$$\{(b/a,c/a)+(d,e): a\in A, (b,c) \in (A \times A) \cap l_{i}, (d,e) \in (A \times A) \cap l_{i+1} \}.$$
There are at least $|A|$ choices for $(b/a,c/a)$ and at least $|A|^2/k$ choices for $(d,e)$. Since all of these sums are distinct, we obtain at least $|A|^3/k$ elements of $(A/A+A) \times (A/A+A)$ lying in between $l_i$ and $l_{i+1}$. Summing over all $1\leq i \leq k-1$, it follows that
$$|A/A+A|^2 \gg |A|^3.$$

There are two rather crude steps in this argument. The first is the observation that there are at least $|A|$ choices for the point $(b/a,c/a)$. In fact, the number of points of this form is equal to the cardinality of product set of $A$ and a set of size $|A|^2/k$. This could be as small as $|A|$, but one would typically expect it to be considerably larger. This extra information was used by Shkredov \cite{shkredov} in his proof of \eqref{a+a:a}.

The second wasteful step comes at the end of the argument, when we only consider sums coming from pairs of lines which are neighbours. This means that we consider only $k-1$ pairs of lines out of a total of ${k \choose 2}$. A crucial ingredient in the proof of \eqref{KS} was the ability to find a way to count sums coming from more than just neighbouring lines.

The proof of Theorem \ref{thm:main} deals with these two steps more efficiently. Ideas from \cite{shkredov} are used to improve upon the first step, and then ideas from \cite{KS} improve upon the second step. We also make use of the fact that the set $A/A$ is invariant under the function $f(x)=1/x$, which allows us to use results on convexity and sumset of Elekes, Nathanson and Ruzsa \cite{ENR} in order to get a better exponent in Theorem \ref{thm:main}.

\section{Notation and Preliminary results}

Throughout the paper, the standard notation
$\ll,\gg$ is applied to positive quantities in the usual way. Saying $X\gg Y$ means that $X\geq cY$, for some absolute constant $c>0$.

The main tool is the Szemer\'{e}di-Trotter Theorem.

\begin{theorem} \label{thm:SzT}
Let $P$ be a finite set of points in $\mathbb R^2$ and let $L$ be a finite set of lines. Then
$$|\{(p,l)\in P \times L : p \in l\}| \ll (|P||L|)^{2/3}+|P|+|L|.$$

\end{theorem}

Define
\begin{equation}
d(A)=\min_{C \neq \emptyset} \frac{|AC|^2}{|A||C|}.
\label{eq:d(A)def}
\end{equation}
We will need the following consequence of the Szemer\'{e}di-Trotter Theorem, which is \cite[Corollary 8]{KS}.
\begin{lemma} \label{thm:ST}
Let $A_1,A_2$ and $A_3$ be finite sets of real numbers and let $\alpha_1,\alpha_2$ and $\alpha_3$ be arbitrary non-zero real numbers. Then the number of solutions to the equation
$$\alpha_1a_1+\alpha_2a_2+\alpha_3a_3=0,$$
such that $a_1 \in A_1$, $a_2 \in A_2$ and $a_3 \in A_3$, is at most 
$$C\cdot d^{1/3}(A_1)|A_1|^{1/3}|A_2|^{2/3}|A_3|^{2/3},$$ for some absolute constant $C$.
\end{lemma}

Another application of (a variant of) the Szemer\'{e}di-Trotter Theorem is the following result of Elekes, Nathanson and Ruzsa \cite{ENR}:

\begin{theorem} \label{ENR} Let $f: \mathbb R \rightarrow \mathbb R$ be a strictly convex or concave function and let $X,Y,Z \subset \mathbb R$ be finite. Then
$$|f(X)+Y||X+Z| \gg |X|^{3/2}|Y|^{1/2}|Z|^{1/2}.$$
\end{theorem}

In particular, this theorem can be applied with $f(x)=1/x$, $X=A/A$, $Y=Z=A$, using the fact that $f(A/A)=A/A$, to obtain the following corollary:
\begin{cor} \label{ENRcor} For any finite set $A \subset \mathbb R$,
$$|A/A+A| \gg |A/A|^{3/4}|A|^{1/2}.$$
\end{cor}

\section{Proof of main theorem}

Recall that the aim is to prove the inequality
$$|A/A+A| \gg \frac{|A|^{\frac{3}{2}+\frac{1}{26}}}{\log^{1/2}|A|}.$$

Consider the point set $A \times A$ in the plane. At the outset, we perform a dyadic decomposition, and then apply the pigeonhole principle, in order to find a large subset of $A\times A$ consisting of points lying on lines through the origin which contain between $\tau$ and $2\tau$ points, where $\tau$ is some real number.

Following the notation of \cite{KS}, for a real number $\lambda$, define
$$ \mathcal A_{\lambda}:= \left\{(x,y) \in A \times A : \frac{y}{x}=\lambda \right\},$$
and its projection onto the horizontal axis,
$$A_{\lambda}:=\{x:(x,y) \in \mathcal A_{\lambda}\}.$$
Note that $|A_{\la}|=|A \cap \la A|$ and
\begin{equation}
\sum_{\la} |A_{\la}|=|A|^2.
\label{obvious}
\end{equation}

Let $S_{\tau}$ be defined by
$$|S_{\tau}|:=|\{\lambda: \tau \leq |A \cap \lambda A| < 2\tau \}|.$$
After dyadically decomposing the sum \eqref{obvious}, we have
$$|A|^2=\sum_{\la} |A_{\la}| = \sum_{j=1}^{\lceil \log|A| \rceil} \sum_{\la \in S_{2^{j-1}}}|A_{\la}| .$$
Applying the pigeonhole principle, we deduce that there is some $\tau$ such that
\begin{equation}
\sum_{\la \in S_{\tau}}|A_{\la}| \geq \frac{|A|^2}{\lceil \log|A| \rceil} \geq \frac{|A|^2}{2\log|A|}.
\label{sumbound}
\end{equation}
Since $\tau \leq |A|$, this implies that
\begin{equation}
|S_{\tau}| \geq \frac{|A|}{2\log |A|}.
\label{Sbound}
\end{equation}
Also, since $|A_{\la}| < 2\tau$ for any $\la \in S_{\tau}$, we have
\begin{equation}
\tau|S_{\tau}| \gg \frac{|A|^2}{ \log|A|}.
\label{taubound}
\end{equation}

\subsection{A lower bound for $\tau$}

Suppose that $|A/A| \geq |A|^{\frac{4}{3}+\frac{2}{39}}$. Then, by Corollary \ref{ENRcor},
$$|A/A+A| \gg |A/A|^{\frac{3}{4}}|A|^{\frac{1}{2}} \gg |A|^{\frac{3}{2}+\frac{1}{26}},$$
as required. Therefore, we may assume that $|A/A| \leq |A|^{\frac{4}{3}+\frac{2}{39}}$. In particular, by \eqref{taubound},
$$\tau|A|^{\frac{4}{3}+\frac{2}{39}} \geq \tau|A/A| \geq \tau|S_{\tau}| \gg \frac{|A|^2}{ \log|A|}.$$
Therefore
\begin{equation}
\tau \gg \frac{|A|^{\frac{2}{3}-\frac{2}{39}}}{ \log|A|}.
\label{taubound3}
\end{equation}

\subsection{An upper bound for $d(A)$}

Define $P$ to be the subset of $A \times A$ lying on the union of the lines through the origin containing between $\tau$ and $2\tau$ points. That is, $P = \cup_{\la \in S_{\tau}} \mathcal A_{\la}$. We will study vector sums coming from this point set by two different methods, and then compare the bounds in order to prove the theorem. To begin with, we use the methods from the paper \cite{shkredov} to obtain an upper bound for $d(A)$. The deduction of the forthcoming bound \eqref{part1} is a minor variation of the first part of the proof of \cite[Theorem 13]{shkredov}.

After carrying out the aforementioned pigeonholing argument, we have a set of $|S_{\tau}|$ lines through the origin, each containing approximately $\tau$ points from $A \times A$. Label the lines $l_1,l_2,\dots,l_{|S_{\tau}|}$ in increasing order of steepness. The line $l_i$ has equation $y=q_ix$ and so $q_1<q_2<\dots<q_{|S_{\tau}|}$. For any $1 \leq i \leq |S_{\tau}|-1$, consider the sum set
\begin{equation}
\mathcal A_{q_i}+\mathcal A_{q_{i+1}} \cdot \Delta(A^{-1}) \subset (A+A/A) \times (A+A/A),
\label{vectorsum}
\end{equation}
where $\Delta(B)=\{(b,b):b \in B\}$. Note that $\mathcal A_{q_{i+1}} \cdot \Delta(A^{-1})$ has cardinality $|A_{q_{i+1}}A^{-1}|$, and therefore the set in \eqref{vectorsum} has at least $|A_{q_{i+1}}A^{-1}||A_{q_i}|$ elements, all of which lie in between $l_i$ and $l_{i+1}$. This is a consequence of the observation of Solymosi that the sum set of $m$ points on one line through the origin and $n$ points on another line through the origin consists of $mn$ points lying in between the two lines. It is important to note that this fact is dependent on the points lying inside the positive quadrant of the plane, which is why the assumption that $A$ consists of strictly positive reals is needed for this proof.

Summing over all $1 \leq i < |S_{\tau}|$, applying the definition of $d(A)$ and using the bounds \eqref{taubound} and \eqref{sumbound}, we obtain
\begin{align*}
|A/A+A|^2 &\geq \sum_{i=1}^{|S_{\tau}|-1} |A_{q_i}||A/A_{q_{i+1}}|
\\&\geq |A|^{1/2}d^{1/2}(A)\sum_{i=1}^{|S_{\tau}|-1}|A_{q_i}||A_{q_{i+1}}|^{1/2}
\\& \gg \frac{|A|^{3/2}d^{1/2}(A)}{|S_{\tau}|^{1/2}\log^{1/2}|A|} \sum_{i=1}^{|S_{\tau}|-1}|A_{q_i}|
\\& \gg \frac{|A|^{7/2}d^{1/2}(A)}{|S_{\tau}|^{1/2}\log^{3/2}|A|}.
\end{align*}
This can be rearranged to obtain
\begin{equation}
d(A) \ll \frac{|A/A+A|^4|S_{\tau}|\log^3|A|}{|A|^7}.
\label{part1}
\end{equation}

This bound will be utilised later in the proof. We now analyse the vector sums in a different way, based on the approach of \cite{KS}.

\subsection{Clustering setup}

For each $\lambda \in S_{\tau}$, we identify an element from $\mathcal A_{\la}$, which we label $(a_{\lambda},\lambda a_{\lambda})$. These fixed points will have to be chosen with a little care later, but for the next part of the argument, we can think of the choice of $(a_{\lambda},\lambda a_{\lambda})$ as completely arbitrary, since the required bound holds whichever choice we make for these fixed points. 

Then, fixing two distinct slopes $\lambda$ and $\lambda'$ from $S_{\tau}$ and following the observation of Balog \cite{balog}, we note that at least $\tau|A|$ distinct elements of $(A/A+A) \times (A/A+A)$ are obtained by summing points from the two lines. Indeed,
$$\mathcal A_{\lambda}+(a_{\lambda'},\lambda'a_{\lambda'}) \cdot \Delta(A^{-1}) \subset (A/A+A) \times (A/A+A).$$ 
Once again, these vector sums are all distinct and have slope in between $\lambda$ and $\lambda'$.

Following the strategy of Konyagin and Shkredov \cite{KS}, we split the family of $|S_{\tau}|$ slopes into clusters of $2M$ consecutive slopes, where $2\leq 2M \leq |S_{\tau}|$ and $M$ is a parameter to be specified later. For example, the first cluster is $U_1= \{l_1,\dots,l_{2M}\}$, the second is $U_2=\{l_{2M+1},\dots,l_{4M}\}$, and so on. We then split each cluster arbitrarily into two disjoint subclusters of size $M$. For example, we have $U_1=V_1 \sqcup W_1$ where $V_1=\{l_1,\dots,l_M\}$ and $W_1=\{l_{M+1},\dots,l_{2M}\}$.

The idea is to show that each cluster determines many different elements of $(A+A/A) \times (A+A/A)$. Since the slopes of these elements are in between the maximal and minimal values in that cluster, we can then sum over all clusters without overcounting.

If a cluster contains exactly $2M$ lines, then it is called a \textit{full cluster}. Note that there are $\left\lfloor \frac{|S_{\tau}|}{2M} \right\rfloor \geq \frac{|S_{\tau}|}{4M}$ full clusters, since we place exactly $2M$ lines in each cluster, with the possible exception of the last cluster which contains at most $2M$ lines.

The proceeding analysis will work in exactly the same way for any full cluster, and so for simplicity of notation we deal only with the first cluster $U_1$. We further simplify this by writing $U_1=U$, $V_1=V$ and $W_1=W$.

Let $\mu$ denote the number of elements of $(A/A+A) \times (A/A+A)$ which lie in between $l_1$ and $l_{2M}$. Then\footnote{For the sake of simplicity of presentation, a small abuse of notation is made here. The lines in $V$ and $W$ are identified with their slopes. In this way, the notation $\la_i \in V$ is used as a shorthand for $\{ (x,y) : y=\la_i x\} \in V$.}



\begin{equation}\mu \geq \tau|A| M^2 - \sum_{\lambda_1, \lambda_3 \in V ,\lambda_2,\lambda_4 \in W: \{\lambda_1,\lambda_2\} \neq \{\lambda_3,\lambda_4\}} \mathcal E(\lambda_1,\lambda_2,\lambda_3,\lambda_4),
\label{mucount2}
\end{equation}
where
$$\mathcal E(\lambda_1,\lambda_2,\lambda_3,\lambda_4):=|\{z \in (\mathcal A_{\lambda_1}+(a_{\lambda_2},\lambda_2a_{\lambda_2})\cdot\Delta(A^{-1}))\cap (\mathcal A_{\lambda_3}+(a_{\lambda_4},\lambda_4a_{\lambda_4})\cdot\Delta(A^{-1})) \}|.$$
In \eqref{mucount2}, the first term is obtained by counting sums from all pairs of lines in $V \times W$. The second error term covers the overcounting of elements that are counted more than once in the first term. 

The next task is to obtain an upper bound for $\mathcal E(\lambda_1,\lambda_2,\lambda_3,\lambda_4)$ for an arbitrary quadruple $(\la_1,\la_2,\la_3,\la_4)$ which satisfies the aforementioned conditions.

Suppose that 
$$z=(z_1,z_2) \in (\mathcal A_{\lambda_1}+(a_{\lambda_2},\lambda_2a_{\lambda_2})\cdot\Delta(A^{-1}))\cap (\mathcal A_{\lambda_3}+(a_{\lambda_4},\lambda_4a_{\lambda_4})\cdot\Delta(A^{-1})).$$
Then
\begin{align*}
(z_1,z_2) &=(a_1,\lambda_1a_1)+(a_{\lambda_2}a^{-1},\lambda_2a_{\lambda_2}a^{-1})
\\&=(a_3,\lambda_3a_3)+(a_{\lambda_4}b^{-1},\lambda_4a_{\lambda_4}b^{-1}),
\end{align*}
for some $a_1 \in A_{\lambda_1}$, $a_3 \in A_{\lambda_3}$ and $a,b \in A$. Therefore,
\begin{align*}
z_1&=a_1+a_{\lambda_2}a^{-1}=a_3+a_{\lambda_4}b^{-1}
\\z_2&=\lambda_1a_1+\lambda_2a_{\lambda_2}a^{-1}=\lambda_3a_3+\lambda_4a_{\lambda_4}b^{-1}
\end{align*}

\subsection{Bounding $\mathcal E(\lambda_1,\lambda_2,\lambda_3,\lambda_4)$ in the case when $\la_4 \neq \la_2$}
Let us assume first that $\la_4 \neq \la_2$. Note that this assumption implies that $\la_4\neq \la_1,\la_2,\la_3$. We have
$$0=\lambda_1a_1+\lambda_2a_{\lambda_2}a^{-1}-\lambda_3a_3-\lambda_4a_{\lambda_4}b^{-1} - \lambda_4(a_1+a_{\lambda_2}a^{-1}-a_3-a_{\lambda_4}b^{-1}),$$
and thus
\begin{equation}
0=a_{\lambda_2}(\lambda_2-\lambda_4)a^{-1}+(\lambda_1-\lambda_4)a_1+(\lambda_4-\lambda_3)a_3.
\label{STsetup}
\end{equation}

Note that the values $\lambda_1-\lambda_4, a_{\lambda_2}(\lambda_2-\lambda_4)$ and $\lambda_4-\lambda_3$ are all non-zero. We have shown that each contribution to $\mathcal E (\lambda_1,\lambda_2,\lambda_3,\lambda_4)$ determines a solution to \eqref{STsetup} with $(a,a_1,a_3) \in A \times A_{\la_1} \times A_{\la_3}$. Furthermore, the solution to \eqref{STsetup} that we obtain via this deduction is unique, and so a bound for  $\mathcal E(\lambda_1,\lambda_2,\lambda_3,\lambda_4)$ will follow from a bound to the number of solutions to \eqref{STsetup}. 

It therefore follows from an application of Lemma \ref{thm:ST} that
$$\mathcal E(\lambda_1,\lambda_2,\lambda_3,\lambda_4) \leq C\cdot d(A^{-1})^{1/3}|A|^{1/3}(A)\tau^{4/3}=C\cdot d(A)^{1/3}|A|^{1/3}(A)\tau^{4/3},$$
where $C$ is an absolute constant. Therefore,
\begin{equation}
\mu \geq M^2 |A|\tau - M^4Cd^{1/3}(A)|A|^{1/3}\tau^{4/3} -\sum_{\lambda_1, \lambda_3 \in V ,\lambda_2 \in W: \lambda_1 \neq \lambda_3} \mathcal E(\lambda_1,\lambda_2,\lambda_3,\lambda_2)
\label{mu}
\end{equation}
We now impose a condition on the parameter $M$ (recall that we will choose an optimal value of $M$ at the conclusion of the proof) to ensure that the first error term is dominated by the main term. We need
$$CM^4d^{1/3}(A)|A|^{1/3}\tau^{4/3} \leq \frac{M^2|A|\tau}{2},$$
which simplifies to 
\begin{equation} \label{Mcond}
M\leq  \frac{|A|^{1/3}}{\sqrt{2C}d^{1/6}(A)\tau^{1/6}} .
\end{equation}
With this restriction on $M$, we now have
\begin{equation}
\mu \geq \frac{M^2 |A|\tau}{2}  -\sum_{\lambda_1, \lambda_3 \in V ,\lambda_2 \in W: \lambda_1 \neq \lambda_3} \mathcal E(\lambda_1,\lambda_2,\lambda_3,\lambda_2).
\label{mu}
\end{equation}
It remains to bound this second error term.

\subsection{Bounding $\mathcal E(\lambda_1,\lambda_2,\lambda_3,\lambda_4)$ in the case $\la_4 = \la_2$}


It is in this case that we need to take care to make good choices for the fixed points $(a_{\lambda},\lambda a_{\la})$ on each line $l_{\la}$.

Fix $\la_2 \in W $. We want to prove that there is a choice for $(a_{\lambda_2},\lambda a_{\la_2}) \in \mathcal A_{\la_2}$ such that
$$\sum_{\lambda_1, \lambda_3 \in V : \lambda_1 \neq \lambda_3}  \mathcal E(\lambda_1,\lambda_2,\lambda_3,\lambda_2) \ll M^2\tau^{1/3}|A|^{4/3}.$$
We will do this using the Szemer\'{e}di-Trotter Theorem. Consider the sum
$$\sum_{a_{\la_2} \in A_{\la_2}} \sum_{\lambda_1, \lambda_3 \in V : \lambda_1 \neq \lambda_3}  \mathcal |\{z \in (\mathcal A_{\lambda_1}+(a_{\lambda_2},\lambda_2a_{\lambda_2})\cdot\Delta(A^{-1}))\cap (\mathcal A_{\lambda_3}+(a_{\lambda_2},\lambda_2a_{\lambda_2})\cdot\Delta(A^{-1})) \}|.$$

Suppose that 
$$z=(z_1,z_2) \in (\mathcal A_{\lambda_1}+(a_{\lambda_2},\lambda_2a_{\lambda_2})\cdot\Delta(A^{-1}))\cap (\mathcal A_{\lambda_3}+(a_{\lambda_2},\lambda_2a_{\lambda_2})\cdot\Delta(A^{-1})).$$
Then
\begin{align*}
(z_1,z_2) &=(a_1,\lambda_1a_1)+(a_{\lambda_2}a^{-1},\lambda_2a_{\lambda_2}a^{-1})
\\&=(a_3,\lambda_3a_3)+(a_{\lambda_2}b^{-1},\lambda_2a_{\lambda_2}b^{-1}),
\end{align*}
for some $a_1 \in A_{\lambda_1}$, $a_3 \in A_{\lambda_3}$ and $a,b \in A$. Therefore,
\begin{align*}
z_1&=a_1+a_{\lambda_2}a^{-1}=a_3+a_{\lambda_2}b^{-1}
\\z_2&=\lambda_1a_1+\lambda_2a_{\lambda_2}a^{-1}=\lambda_3a_3+\lambda_2a_{\lambda_2}b^{-1}.
\end{align*}

We have
$$0=\lambda_1a_1+\lambda_2a_{\lambda_2}a^{-1}-\lambda_3a_3-\lambda_2a_{\lambda_2}b^{-1} - \lambda_1(a_1+a_{\lambda_2}a^{-1}-a_3-a_{\lambda_2}b^{-1}),$$
and thus
\begin{equation}
\frac{\la_3-\la_1}{\la_2-\la_1}a_3=a_{\lambda_2}(a^{-1}-b^{-1}).
\label{STsetup2}
\end{equation}

As in the previous subsection, this shows that the quantity 
$$\sum_{a_{\la_2} \in A_{\la_2}}  \sum_{\la_1 \neq \la_3 \in V }|(\mathcal A_{\lambda_1}+(a_{\lambda_2},\lambda_2a_{\lambda_2})\cdot\Delta(A^{-1}))\cap (\mathcal A_{\lambda_3}+(a_{\lambda_2},\lambda_2a_{\lambda_2})\cdot\Delta(A^{-1}))|$$ 
is no greater than the number of solutions to \eqref{STsetup2} such that $(\la_1,\la_3,a,b,a_{\lambda_2},a_3) \in V \times V \times A \times A \times A_{\la_2} \times A_{\la_3}$.

Fix, $\la_1, \la_3 \in V$ such that $\la_1\neq \la_3$. Let $Q=A^{-1} \times A_{\la_3}$. Define $l_{m,c}$ to be the line with equation $\frac{\la_3-\la_1}{\la_2-\la_1}y=m(x-c)$ and define $L$ to be the set of lines
$$L=\{l_{a_{\la_2},b^{-1}}: a_{\la_2} \in A_{\la_2}, b \in A \}.$$
Note that $|Q|\approx |L| \approx \tau|A|$ and so
$$I(Q,L) \ll (\tau|A|)^{4/3}.$$
Repeating this analysis via the Szemer\'{e}di-Trotter Theorem for each pair of distinct $\la_1,\la_3 \in V$, it follows that the number of solutions to \eqref{STsetup2} is $O(M^2(\tau|A|)^{4/3})$. In summary,
$$\sum_{a_{\la_2} \in A_{\la_2}} \sum_{\la_1 \neq \la_3 \in V }|(\mathcal A_{\lambda_1}+(a_{\lambda_2},\lambda_2a_{\lambda_2})\cdot\Delta(A^{-1}))\cap (\mathcal A_{\lambda_3}+(a_{\lambda_2},\lambda_2a_{\lambda_2})\cdot\Delta(A^{-1}))|\ll M^2(\tau|A|)^{4/3}.$$
Therefore, by the pigeonhole principle, there is some $a_{\la_2} \in A_{\la_2}$ such that
\begin{equation} \label{fix1}
 \sum_{\la_1 \neq \la_3 \in V }|(\mathcal A_{\lambda_1}+(a_{\lambda_2},\lambda_2a_{\lambda_2})\cdot\Delta(A^{-1}))\cap (\mathcal A_{\lambda_3}+(a_{\lambda_2},\lambda_2a_{\lambda_2})\cdot\Delta(A^{-1}))|\ll M^2\tau^{1/3}|A|^{4/3}.
 \end{equation}
We can then choose the fixed point $(a_{\la_2}, \la_2a_{\la_2})$ on $l_{\la_2}$ to be that corresponding to the value $a_{\la_2}$ satisfying inequality \eqref{fix1}. This in fact shows that
\begin{equation} \label{fix2}
 \sum_{\la_1 \neq \la_3 \in V }\mathcal E (\la_1,\la_2,\la_3,\la_2) \ll M^2\tau^{1/3}|A|^{4/3}.
 \end{equation}
We repeat this process for each $\la_2 \in W$ to choose a fix point for each line with slope in $W$. Summing over all $\la_2 \in W$, we now have
\begin{equation} \label{fix2}
 \sum_{\la_1 \neq \la_3 \in V, \la_2 \in W }\mathcal E (\la_1,\la_2,\la_3,\la_2) \ll M^3\tau^{1/3}|A|^{4/3}.
 \end{equation}
 
We have a bound for error term in \eqref{mu}. Still, we need to impose a condition on $M$ so that this error term is dominated by the main term. We need
$$M^3\tau^{1/3}|A|^{4/3} \leq \frac{M^2|A|\tau}{4},$$
which simplifies to 
\begin{equation} \label{Mcond2}
M\leq  \frac{\tau^{2/3}}{4|A|^{1/3}} .
\end{equation}
With this restriction on $M$, we now have
\begin{equation}
\mu \geq \frac{M^2|A|\tau}{4}.
\label{almost}
\end{equation}

Our integer parameter $M$ must satisfy \eqref{Mcond} and \eqref{Mcond2}. We therefore choose
$$M:=\min \left\{ \left \lfloor \frac{|A|^{1/3}}{\sqrt{2C}d^{1/6}(A)\tau^{1/6}} \right \rfloor, \left \lfloor \frac{\tau^{2/3}}{4|A|^{1/3}} \right \rfloor \right\}.$$
Summing over the full clusters, of which there are at least $\frac{|S_{\tau}|}{4M}$, yields
\begin{align}
|A/A+A|^2&\geq \frac{|S_{\tau}|}{4M}\frac{M^2}{4}|A|\tau 
\\&\gg |S_{\tau}|M|A|\tau 
\label{aa+a} 
\end{align}

\subsection{Choosing $M$ - case 1}

Suppose first that $M=\left \lfloor \frac{|A|^{1/3}}{\sqrt{2C}d^{1/6}(A)\tau^{1/6}} \right \rfloor$.

Recall that we need $2 \leq 2M \leq |S_{\tau}|$. It is easy to check that the upper bound for $M$ is satisfied. Indeed, it follows from \eqref{Sbound} that
$$2M \leq \frac{2}{\sqrt{2C}}|A|^{1/3} \leq \frac{|A|}{2 \log |A|} \leq  |S_{\tau}|,$$
The first inequality above uses the fact that $d(A) \geq 1$ for all $A$ (since one can take $C$ to be a singleton in the \eqref{eq:d(A)def}), as well as the bound $\tau \geq 1$.
The second inequality is true for sufficiently large $|A|$. Since smaller sets can be dealt with by choosing sufficiently small implied constants in the statement, we may assume that $2M \leq |S_{\tau}|$.

Assume first that $M \geq 1$ (we will deal with the other case later). Then, by \eqref{aa+a} and the definition of $M$
$$|A/A+A|^2 \gg \frac{|S_{\tau}||A|^{4/3}\tau^{5/6}}{d^{1/6}(A)}.$$
Applying the inequality $|S_{\tau}|\tau \gg \frac{|A|^2}{\log |A|}$, it follows that
\begin{equation}
d^{1/6}(A)|A/A+A|^2 \gg \frac{ |A|^{3}|S_{\tau}|^{1/6}}{\log^{5/6}|A|}.
\label{aa+a2}
\end{equation}
After bounding the left hand side of this inequality using \eqref{part1}, we obtain
$$\frac{|A/A+A|^{2/3}|S_{\tau}|^{1/6}\log^{1/2}|A|}{|A|^{7/6}}|A/A+A|^2 \gg d^{1/6}(A)|A/A+A|^2 \gg \frac{ |A|^{3}|S_{\tau}|^{1/6}}{\log^{5/6}|A|}.$$
Rearranging this expression leads to the bound
$$|A/A+A| \gg \frac{|A|^{25/16}}{\log^{1/2}|A|},$$
which is stronger than the claim of the theorem.


It remains is to consider what happens if $M \leq 1$. Indeed, if this is the case, then
$$\frac{|A|^{1/3}}{\sqrt{8C}d^{1/6}(A)\tau^{1/6}}<1$$
and so
$$\frac{|A|^{1/3}}{d^{1/6}(A)\tau^{1/6}} \ll 1.$$
After applying the bound $\tau \leq |A|^2/|S_{\tau}|$, it follows that
$$\frac{|S_{\tau}|^{1/6}}{d^{1/6}(A)} \ll 1 \ll \frac{|A/A+A|^2}{|A|^3},$$
where the latter inequality is a consequence of Theorem \ref{thm:old}. In particular, this implies that \eqref{aa+a2} holds. We can then repeat the earlier analysis and once again reach the conclusion that 
$$|A/A+A| \gg \frac{|A|^{\frac{3}{2}+\frac{1}{16}}}{\log^{1/2}|A|}.$$

\subsection{Choosing $M$ - case 2}

Suppose now that $M=\left \lfloor \frac{\tau^{2/3}}{4|A|^{1/3}} \right \rfloor$.

Again, we need to check that $2 \leq 2M \leq |S_{\tau}|$. If the lower bound does not hold then \eqref{taubound3} gives a contradiction for sufficiently large $|A|$. Smaller sets can be dealt with by choosing sufficiently small implied constants in the statement. If the upper bound does not hold then
$$\frac{\tau^{2/3}}{|A|^{1/3}} \geq 2M>|S_{\tau}|.$$
Multiplying both sides of this inequality by $\tau$ and applying \eqref{taubound} gives the contradiction
$$|A|^{5/3} \geq \tau^{5/3} \gg \frac{|A|^{7/3}}{\log |A|}.$$


Since this choice of $M$ is valid, we can now conclude the proof. From \eqref{taubound3}, we have
$$M \gg \frac{\tau^{2/3}}{|A|^{1/3}} \gg \frac{|A|^{\frac{1}{13}}}{\log^{\frac{2}{3}}|A|}.$$
Then, by \eqref{aa+a} and \eqref{taubound},
$$|A/A+A|^2 \gg \frac{|A|^3}{\log|A|}M \gg \frac{|A|^{3+\frac{1}{13}}}{\log^{5/3}|A|} .$$
We conclude that
$$|A/A+A| \gg \frac{|A|^{\frac{3}{2}+\frac{1}{26}}}{\log^{5/6}|A|},$$
and so the proof is complete. \qedsymbol


\begin{thebibliography}{99}

\bibitem{balog} A. Balog `A note on sum-product estimates', \textit{Publ. Math. Debrecen} \textbf{79}, no. 3-4 (2011), 283-289.








	

\bibitem{ENR} G. Elekes, M. Nathanson and I. Ruzsa, `Convexity and sumsets', \textit{J. Number Theory.} \textbf{83} (1999), 194-201.






\bibitem{KS} S. Konyagin and I. Shkredov, `On sum sets of sets, having small product set', \textit{Proc. Steklov Inst. Math.} \textbf{290} (2015), 288-299.

\bibitem{KS2} S. Konyagin and I. Shkredov, `New results on sums and products in $\mathbb R$', \textit{Proc. Steklov Inst. Math.} \textbf{294} (2016), 87-98.




\bibitem{MORNS} B. Murphy, O. Roche-Newton and I. Shkredov `Variations on the sum-product problem', \textit{SIAM J. Discrete Math.} 29 (2015), no. 1, 514-540.




\bibitem{RRSS} O. Roche-Newton, I. Z. Ruzsa, C.-Y. Shen and I. D. Shkredov `On the size of the set $AA+A$', \textit{To appear in JLMS}.


\bibitem{RSS} M. Rudnev,  I. D. Shkredov and S. Stevens, 
`On the energy variant of the sum-product conjecture', \textit{arxiv:1607.05053} (2016). 












\bibitem{shkredov} I. Shkredov, `On a question of A. Balog', \textit{Pacific J. Math.} 280 (2016), no. 1, 227-240.

\bibitem{solymosi} J. Solymosi, `Bounding multiplicative energy by the sumset', \textit{Adv. Math.} \textbf{222} (2009), 402-408.


\bibitem{tv} T. Tao, V. Vu. 'Additive combinatorics' \textit{Cambridge University Press} (2006).









\end{thebibliography}
\end{document}